\documentclass[preprint,12pt]{elsarticle}

\usepackage{graphicx}
\usepackage{amssymb}
\usepackage[noend]{algorithmic}
\usepackage{algorithm}

\newtheorem{thm}{Theorem}
\newtheorem{lem}[thm]{Lemma}
\newtheorem{props}[thm]{Proposition}
\newtheorem{cors}[thm]{Corollary}
\newdefinition{rmk}{Remark}
\newdefinition{defs}[rmk]{Definition}
\newproof{pf}{Proof}
\renewcommand{\algorithmicrequire}{\textbf{Input:}}
\renewcommand{\algorithmicensure}{\textbf{Output:}}

\journal{arXiv}

\begin{document}

\begin{frontmatter}

\title{The bottleneck $2$-connected $k$-Steiner network problem for $k\leq 2$}

\author[]{M.~Brazil}

\author[]{C.J.~Ras\corref{cor2}}

\author[]{D.A.~Thomas}

\begin{abstract}
The geometric bottleneck Steiner network problem on a set of vertices $X$ embedded in a normed plane requires one to construct a graph $G$
spanning $X$ and a variable set of $k\geq 0$ additional points, such that the length of the longest edge is minimised. If no other constraints
are placed on $G$ then a solution always exists which is a tree. In this paper we consider the Euclidean bottleneck Steiner network problem for
$k\leq 2$, where $G$ is constrained to be $2$-connected. By taking advantage of relative neighbourhood graphs, Voronoi diagrams, and the tree
structure of block cut-vertex decompositions of graphs, we produce exact algorithms of complexity $O(n^2)$ and $O(n^2\log n)$ for the cases
$k=1$ and $k=2$ respectively. Our algorithms can also be extended to other norms such as the $L_p$ planes.

\end{abstract}


\end{frontmatter}

\section{Introduction}
In communication networks a \textit{bottleneck} can be any node or link at which a performance objective attains its least desirable value. For
instance, in wireless sensor networks we may define a bottleneck parameter on the network as the length of the longest edge (link), where the
benefit of minimising the length of a link comes from the observation that the energy consumption of the incident transmitting node, for each
transmission, increases with the length of the link. Due to the requirement of prolonged autonomy in wireless sensor networks, and the
subsequent use of batteries, optimisation of power in individual nodes is a primary goal. This particular bottleneck parameter is therefore a
common optimisation objective in the modelling of sensor network deployments. Graph models dealing with the minimisation of the longest edge
also have wide applicability in other areas, for instance in VLSI layout, general communication network design, and location problems; see
\cite{sarr} for an introduction to this topic.

Previous work on the longest edge minimisation problem in graphs has centred on properties and algorithms for the construction of
\textit{bottleneck Steiner trees}, both in the geometric version of the problem, and in the graph version where solutions are required to be
subgraphs of a given weighted graph. In all versions of the problem one is required to construct a spanning tree on a given set of $n$ vertices
such that the longest edge has minimum length (or weight), and where a set of additional points (called Steiner points) are available during the
construction. In geometric versions Steiner points can generally be located anywhere in the plane, and therefore, to ensure that the bottleneck
can not be made arbitrarily small, an upper bound $k$ is placed on their total number. In the Euclidean and rectilinear planes, and also in
graphs, the problem has been shown to be NP-hard; see \cite{berm,sarr,wang}. Recent papers provide exact algorithms for the $L_p$ metric and
other normed planes; for instance \cite{bae1,bae2,brazil}. In particular, in \cite{bae1} Bae et al. present an $O(f(k)\cdot (n^k+n\log n))$
algorithm for $L_p$ metrics with $1<p<\infty$, where $n$ is the number of non-Steiner vertices and $f(k)$ is a function of $k$ only. They make
use of a technique based on smallest colour spanning disks and farthest colour Voronoi diagrams, which we also employ for our algorithms.

As a model for wireless network deployment the bottleneck Steiner tree problem is only an initial step towards the more general (and realistic)
aim of modelling networks of higher connectivity. The benefits of multi-path connectivity in networks are numerous, and include robustness and
survivability of the network in the event of node failure. In wireless sensor networks another benefit of multiple available paths is the
possibility of diverting traffic when a node's available power is low, and the subsequent extension of the lifetime (or time till first
maintenance) of the network.

Few results exist in the literature for the bottleneck Steiner problem when the solution graph is required to be anything other than a tree. The
case when the resultant graph is required to be $2$-connected, but no Steiner points are allowed, finds application as a heuristic for the
bottleneck Travelling Salesman Problem, as was shown by Timofeev in \cite{timo} and by Parker and Rardin in \cite{park}. Various authors (see
\cite{mschang,manku,punn}) consequently produced fast polynomial algorithms for the so called \textit{bottleneck biconnected spanning subgraph
problem}, the fastest of which provides an $O(m)$ exact algorithm when the initial given graph contains $m$ edges. This translates into an
$O(n^2)$ algorithm for the geometric problem, where all edges of the complete graph are assumed to be available.

This paper presents algorithms for solving the bottleneck Steiner problem in the Euclidean plane when the solution graph is required to be
$2$-connected and contains exactly $k=1$ or $k=2$ Steiner points. We discover new properties of bottleneck Steiner $2$-connected networks that
are based on the well-known block cut-vertex decomposition of graphs. This allows us to develop an $O(n^2$) algorithm for solving the problem
when $k=1$, and an $O(n^2\log n)$ algorithm when $k=2$. We also provide an outline of the generalisation of our techniques to other planar
norms.

The paper is divided into three main parts. Section \ref{sec1} deals with notation and provides a few structural results that are relevant to
both cases $k=1,2$. In Section \ref{sec2} we focus on the case $k=1$ and in Section \ref{sec3} on the case $k=2$.

\section{Notation \& Preliminaries}\label{sec1}
Throughout this paper we only consider finite, simple, and undirected graphs. Let $X$ be a set of vertices embedded in $\mathbb{R}^2$. If
$G=\langle V(G),E(G)\rangle$ is a graph on $X$ then $V(G)=X$ is the \textit{vertex-set} and $E(G)\subset X^2$ the \textit{edge-set} of $G$. If
$A$ is a set of vertices or a graph, and $e$ is an edge incident to some vertex of $A$ then we say $e$ is \textit{incident to} $A$. Two graphs
(or vertex sets) are \textit{adjacent} if there exists an edge incident to both graphs. Two sets of vertices or edges are \textit{independent}
if they are not adjacent or incident to one another. If $G,G^\prime$ are any two graphs, $E\subseteq E(G)$, and $V\subseteq V(G)$, then
$G-E:=\langle V(G),E(G)-E\rangle$, $G-V:=\left\langle V(G)-V,E(G)-\{uv|u\in V\mathrm{\ or\ } v\in V\}\right\rangle$, and $G\cup G^\prime
:=\langle V(G)\cup V(G^\prime), E(G)\cup E(G^\prime)\rangle$.

A graph $G$ is \textit{connected} if there exists a path connecting any pair of vertices in $G$. An \textit{isolated component} is a maximal (by
inclusion) connected subgraph. A \textit{cut-set} $A$ of $G$ is any set of vertices such that $G-A$ has strictly more isolated components than
$G$; if $|A|=1$ then $A$ is a \textit{cut-vertex}. Set $A$ \textit{separates} $W$ from $Z$ in $G$, where $W,Z$ are subgraphs of $G$, if every
path connecting a vertex of $W$ to a vertex of $Z$ contains a vertex of $A$. If $A$ separates any subgraphs of $G$ then $A$ is a cut-set of $G$.

The \textit{vertex-connectivity} or simply \textit{connectivity} $c=c(G)$ of a graph $G$ is the minimum number of vertices whose removal results
in a disconnected or trivial graph. Therefore $c$ is the minimum cardinality of a cut-set of $G$ if $G$ is connected but not complete; $c=0$ if
$G$ is disconnected; and $c=n-1$ if $G=K_n$, where $K_n$ is the complete graph on $n$ vertices. A graph $G$ is said to be
$c^\prime$-\textit{connected} if $c(G)\geq c^\prime$ for some non-negative integer $c^\prime$. In this paper we make an exception for the
connectivity definitions of $K_1,K_2$: we assume that $c(K_1)=c(K_2)=2$. If $G$ is not $K_1$ or $K_2$ then, as a consequence of Menger's
theorem, $G$ is $c^\prime$-connected if and only if for every pair $u,v$ of distinct vertices there are at least $c^\prime$ internally disjoint
$u-v$ paths in $G$. If $G$ is a $2$-connected graph of order at least $3$ then for every triple of vertices of $G$ there exists a cycle
containing them.

A \textit{critical edge} of a $2$-connected graph is an edge such that its removal reduces the graphs connectivity. From \cite{dirac} we know
that an edge is critical if and only if it is not a chord of any cycle. A \textit{block} is a maximal $2$-connected subgraph. The next result is
implicit in many of the proofs in this paper.

\begin{thm}\label{monma1}(see \cite{monma}) Let $G=\langle V,E\rangle$ be a $2$-connected graph with $G^\prime=\langle V^\prime,E^\prime\rangle$ a
subgraph of $G$ induced by $V^\prime$. Then replacing $E^\prime$ in $G$ by any collection of edges $E^{\prime\prime}$ defined on $V^\prime$,
where $G^{\prime\prime}=\langle V^\prime,E^{\prime\prime}\rangle$ is $2$-connected, results in a graph $G^*=\langle V,(E\backslash E^\prime)\cup
E^{\prime\prime})$ which is $2$-connected.
\end{thm}

For any graph $G$ we denote the longest edge of $G$ (where ties have been broken) by $e_{\mathrm{max}}(G)$ and its length by
$\ell_{\mathrm{max}}(G)$.

\begin{defs}
The Euclidean \textit{bottleneck $c$-connected $k$-Steiner network problem} requires one to construct a $c$-connected network $N_k$ spanning $X$
and a set $S_k$ of $k$ Steiner points, such that the $\ell_{\mathrm{max}}(N_k)$ is a minimum across all such networks. The variables are the set
$S_k$ and the topology of the network.
\end{defs}

An optimal solution to the problem is called a \textit{minimum bottleneck $c$-connected $k$-Steiner network}, or $(c,k)$-MBSN. Note that a
$(c,0)$-MBSN is a minimum bottleneck \textit{spanning} $c$-connected network. For the rest of the paper we focus on the case $c=2$ with $k=1,2$.
We also assume throughout that $|X|=n\geq 2$.

Let $\{E_i\}$ be a partition of $E(G)$ into equivalence classes such that two edges are in the same equivalence class if and only if they belong
to a common cycle of $G$. Let $\mathcal{Y}(G)=\{Y_i\}$ where $Y_i$ is the subgraph of $G$ induced by $E_i$. As observed in \cite{hop}, the
partition is well defined; each $Y_i$ is a block of $G$; each non-cut-vertex of $G$ is contained in exactly one of the $Y_i$; each cut-vertex of
$G$ occurs at least twice amongst the $Y_i$; and for each $i,j,i\neq j$, $V(Y_i)\cap V(Y_j)$ consists of at most one vertex, and this vertex (if
it exists) is a cut vertex of $G$. The set $\mathcal{Y}(G)$ is called the \textit{block cut forest} (BCF) of $G$. If $Y_i$ contains exactly one
cut-vertex of $G$ then $Y_i$ is a \textit{leaf block}. An \textit{isolated block} contains no cut-vertices of $G$, i.e., it is a $2$-connected
isolated component of $G$. We use $\mathcal{Y}_0(G)$ to denote the set of leaf blocks of $G$. The \textit{interior} of block $Y_i$, denoted
$Y_i^*$, is the set of all vertices of $Y_i$ that are not cut-vertices of $G$. The unique cut-vertex of $G$ belonging to $Y_i\in
\mathcal{Y}_0(G)$ is denoted by $\tau(Y_i)$.

\begin{thm}(see \cite{tarjan}) The BCF of a graph $G$ with $m$ edges can be constructed in time $O(m)$. As part of the construction
we can calculate the connectivity of $G$, and all leaf blocks as well as all cut-vertices and the blocks that contain them can be specified.
\end{thm}

We define a counter, $b(\cdot)$, as follows. Let $\{G_i\}$ be the set of isolated components of $G$. If $G_i$ is an isolated block then let
$b(G_i)=2$, else let $b(G_i)=|\mathcal{Y}_0(G_i)|$. Finally, let $b(G)=\sum b(G_i)$. Essentially $b(G)$ is the number of leaf blocks plus twice
the number of isolated blocks occurring in $G$ (recall that isolated vertices and isolated edges are blocks according to our definition).

\begin{lem}\label{lembleaf}If $G_1$ is an edge subgraph of $G_2$ then $b(G_1) \geq b(G_2)$.
\end{lem}
\begin{pf}
Every leaf-block of $G_2$ contains a leaf-block or an isolated component of $G_1$. Every isolated block of $G_2$ contains at least two
leaf-blocks or an isolated block of $G_1$.$\ \square$\end{pf}

Let $e$ be any edge of a plane embedded graph. The \textit{lune specified by} $e$ is the region of intersection of the two circles of radius
$|e|$ centred at the endpoints of $e$. Next we define a useful graph for dealing with $2$-connected bottleneck problems.

\begin{defs}(see \cite{mschang}) The $2$-relative neighbourhood graph on $X$ (or $2$-RNG) is the graph $R$ such that $e\in E(R)$ if and only if the
lune specified by $e$ contains (strictly within its boundary) fewer than two vertices of $X$.
\end{defs}

\begin{thm}(see \cite{mschang}) Let $R$ be the $2$-RNG on a given set $X$, with $|X|=n$. Then
\begin{enumerate}
    \item $R$ is 2-connected.
    \item $R$ can be constructed in time $O(n^2)$.
    \item The number of edges of $R$ is $O(n)$.
    \item There exists a $(2,0)$-MBSN, say $N_0$, on $X$ which is a subgraph of $R$. If $R$ is given $N_0$ can be constructed in a time of $O(n\log n)$.
\end{enumerate}
\end{thm}

The algorithms we develop in this paper for constructing $(2,k)$-MBSNs contain a procedure that essentially extends a subgraph $G$ of the
$2$-RNG on $n$ vertices to a $(2,0)$-MBSN containing $G$ as a subgraph and also spanning $k$ variable Steiner points, such that the length of
the longest edge is minimised across all such $(2,0)$-MBSNs. We formalise this concept as follows. Let $G$ be a graph embedded in $\mathbb{R}^2$
and consider the following three variable sets: $S_k=\{s_1,...,s_k\}$, which is a set of $k$ distinct Steiner points in $\mathbb{R}^2$;
$E_S\subset S_k^2$; and $\mathcal{V}=\{V_1,...,V_k\}$, which is a set of subsets of $X=V(G)$. Let $H=\langle V(H),E(H)\rangle$ where $V(H)=X\cup
S_k$ and $E(H)=E(G)\cup E_S\cup\{s_ix_j\ |\ 1\leq i\leq k,x_j\in V_i\}$. If $H$ is $2$-connected then we call $H$ a $k$\textit{-block closure}
of $G$. If $\ell_{\max}(H)\leq \ell_{\max}(H^\prime)$ for any $k$-block closure $H^\prime$ of $G$, then $H$ is an \textit{optimal} $k$-block
closure of $G$. Note that there may be many distinct optimal $k$-block closures for $G$.

A $k$-block closure exists for any graph $G$ when $k\geq 2$: let $S_k$ be any set of $k$ distinct points in the plane, let $E_S=\{(s_i,s_j)\,|\,
i<j\}$, let $V_i=X$ for every $i$, and define $H$ as before. Clearly $H$ is $k$-block closure of $G$. No $1$-block closure exists for a
disconnected graph, since, for any choice of $S_1,E_S$ and $\mathcal{V}$, the resultant $H$ will either be disconnected or the Steiner point
will be a cut-vertex. Observe that $N_k$ is an optimal $k$-block closure of $N_k-S_k$ whenever $N_k$ is a $(2,k)$-MBSN on $X$ with Steiner point
set $S_k$. Therefore $N_1-S_1$ is always connected but $N_2-S_2$ need not be.

A \textit{Steiner edge} is an edge incident to a Steiner point, and for any graph or vertex set $M$ a \textit{Steiner $M$-edge} is an edge
incident to both $S_k$ and $M$. The next lemma is fundamental to our algorithms.

\begin{lem}\label{lemLeaf}For every leaf-block $Y$ of $G$ there exists at least one Steiner $Y^*$-edge in any $k$-block closure of $G$.
\end{lem}
\begin{pf}
If this is not true then $\tau(Y)$ is a cut-vertex of the $k$-block closure, which is a contradiction. $\ \square$\end{pf}

In this paper the construction of an optimal $k$-block closure will usually involve \textit{smallest colour-spanning disks} (SCSDs). Given a
partition of a set $X$ into $\{V_i\}$ where each $V_i$ is assigned a unique colour, an SCSD is a circle of minimum radius that contains at least
one point of each colour. If $|X|=n$ and $|\{V_i\}|$ is constant then an SCSD $C$ can be found in time $O(n\log n)$; see \cite{abel,bae2}.
Clearly $C$ is determined by either two diametrically opposite points, or by three points. These points are referred to (in \cite{bae2}) as the
\textit{determinators} of $C$. The precise way in which one uses SCSDs to construct an optimal $k$-block closure depends on the value of $k$,
and will be discussed in the relevant section.

Proposition~\ref{specN}, below, essentially specifies a useful canonical form for a $(2,k)$-MBSN for any set $X$. The corollaries to this
proposition allow us to greatly minimise the time-complexity of our algorithm for $(2,k)$-MBSN construction later in this section. Before
proving the proposition we require the following lemma.

\begin{lem}\label{edgeswap}Let $N=\langle V,E \rangle$  be a $2$-connected graph. Let $v\in V$ be a vertex of degree $3$ or more in $N$, with
neighbours $x_1$ and $x_2$ such that $vx_1$ and $vx_2$ are critical. Then $x_1x_2 \not\in E$; and replacing $vx_1$ by $x_1x_2$ in $N$ results in
a graph that is also $2$-connected.
\end{lem}
\begin{pf}
Suppose $x_1x_2 \in E$. Since $N$ is $2$-connected and $|V|\geq 4$ it follows that either $vx_1$ or $vx_2$  is a chord of a cycle in $N$,
contradicting the assumption that both edges are critical. Thus, by contradiction, $x_1x_2 \not\in E$.

Let $x_3$ be a third neighbour of $v$ in $N$, other than $x_1$ and $x_2$. Since $N$ is $2$-connected, there exists a path $P_{12}$ between $x_1$
and $x_2$ in $N$ not containing $v$ and there exists a path $P_{23}$ between $x_2$ and $x_3$ in $N$ not containing $v$. The paths $P_{12}$ and
$P_{23}$ are not internally disjoint, since otherwise $vx_2$ would be a chord of the cycle formed by $P_{12},P_{23}, x_3v$ and $vx_1$,
contradicting the assumption that $vx_2$ is critical. It follows that replacing $vx_1$ in $G$ by $x_1x_2$ does not create a cut-vertex at $x_2$.
Clearly no other vertices can become cut-vertices after the replacement, hence the new graph is also $2$-connected. $\ \square$\end{pf}

\begin{figure}[htb]
  \begin{center}
    \includegraphics[scale=0.4]{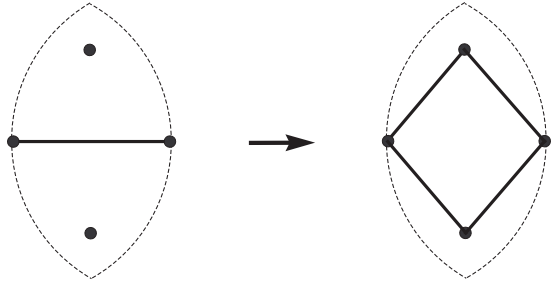}\\
  \end{center}
  \caption{Lune edge-replacement procedure}
  \label{figLune}
\end{figure}

\begin{props}\label{specN}There exists a $(2,k)$-MBSN $N_k$ on $X$, such that $N_k$ is a subgraph of the $2$-RNG on
$V(N_k)$ and the degree of $v$ is at most $5$ for every $v\in V(N_k)$.
\end{props}
\begin{pf}
Let $N$ be any $(2,k)$-MBSN on $X$ such that every edge of $N$ is critical. We also assume that $|V(N_k)| >3$, since otherwise the proposition
is trivially true. The proof is based on running two modification procedures on the edges of $N$, neither of which reduces the connectivity of
the graph: the first reduces the degree of every vertex to at most $5$; the second replaces each edge of $N$ not in the $2$-RNG on $V(N_k)$ by
up to four shorter edges. After each procedure the property of every edge being critical can be maintained by simply deleting any non-critical
edges. We will see that the first procedure does not increase the length of the longest edge in $N$, while, in the second, each edge removed
from $N$ is replaced by shorter edges. It follows that if we alternately run these two modification procedures on $N$, the process must stop
after a finite series of steps, at which point both properties in the proposition have been achieved. It remains to describe the two procedures
and show that each results in a graph that is still $2$-connected.

\noindent \emph{Modification Procedure 1}. Let $v$ be a vertex of $N$ of degree $6$ or more, and let $x_1$ and $x_2$ be two neighbours of $v$
for which $\angle x_1vx_2$ is minimum. We assign the labels to these two neighbours so that $|x_1v| \geq |x_2v|$. Suppose that either $\angle
x_1vx_2 < 60^\circ$ or  $\angle x_1vx_2 = 60^\circ$ and $|x_1v| > |x_2v|$. Then in either case $|x_1x_2| < |x_1v|$, so replacing the edge $x_1v$
by $x_1x_2$ reduces the degree of $v$ and does not increase the length of the longest edge in $N$, but maintains the $2$-connectivity of $N$, by
Lemma~\ref{edgeswap}. Repeating this replacement for every suitable triple $v, x_1, x_2$ results in a graph where a vertex $v$ can only have
degree $6$ if its six neighbours are all equidistant, and each angle between neighbouring pairs of incident edges at $v$ is $60^\circ$. For such
a vertex $v$ we call the subgraph induced by $v$ and its six neighbours a \emph{regular 6-star}. We need to show that we can replace an edge in
$N$ to reduce the degree of the vertex at the centre of a regular 6-star, without creating another regular 6-star elsewhere in the new graph.
Suppose $x_1, x_2$ and $x_3$ are neighbouring vertices in anti-clockwise order to $v$, which is the centre of a regular 6-star, such that
$\angle x_1vx_2 = 60^\circ$, $\angle x_2vx_3 = 60^\circ$ and the edges $vx_i, \ i\in \{1,2,3\}$  are critical. Note that the latter condition
implies that $x_1x_2 \not\in E(N)$ and $x_2x_3 \not\in E(N)$. Suppose we replace $vx_1$ by $x_1x_2$; then, by Lemma~\ref{edgeswap}, the new
graph is still $2$-connected and clearly has the same bottleneck length and total edge length as $N$. But there is no longer a regular 6-star at
$v$, nor has a regular 6-star been created at $x_2$ since $x_2x_3 \not\in E(N)$.

\noindent \emph{Modification Procedure 2}. The second procedure replaces an edge by a 4-cycle if and only if the lune determined by the edge
contains at least $2$ nodes. This procedure replaces the edge by edges of length strictly less than the original edge (see Fig. \ref{figLune}).
The process is described in more detail in \cite{mschang}, where it is also shown that the procedure  maintains the $2$-connectivity of $N$.

Therefore the alternation between these two procedures must eventually terminate and produce a block $N^\prime$ satisfying both conditions. At
this stage we  let $N_k=N^\prime$, completing the proof. $\ \square$\end{pf}

In the rest of this paper we assume that $N_k$ is a $(2,k)$-MBSN on $X$, with Steiner point set $S_k$, satisfying the Proposition~\ref{specN}.
An \textit{external Steiner edge} is a Steiner edge with one end-point not in $S_k$. Let $d$ be the number of external Steiner edges of $N_k$.
For any $G$ we denote the edge-subgraph of $G$ containing all edges of $G$ of length at most $r$ by $G(r)$. Let $R$ be the $2$-RNG on $X$ and
let $\overline{N_k}:=N_k-S_k$. Clearly $\overline{N_k}$ is a subgraph of $R(\ell_{\mathrm{max}}(\overline{N_k}))$.

\begin{cors}\label{maxdeg}$b(R(\ell_{\mathrm{max}}(\overline{N_k})))\leq b(\overline{N_k})\leq d\leq 5k$
\end{cors}
\begin{pf}
The first inequality holds by Lemma \ref{lembleaf} and the second by Lemma \ref{lemLeaf}. The final inequality holds since the degree of any
Steiner point in $N_k$ is at most $5$. $\ \square$\end{pf}

\begin{cors}\label{propR}Let $G=R(\ell_{\mathrm{max}}(\overline{N_k}))$ and let $G^+$ be any optimal $k$-block closure of $G$. Then $G^+$ is
a $(2,k)$-MBSN on $X$.
\end{cors}
\begin{pf}
Since $\overline{N_k}$ is a subgraph of $G$, any $k$-block closure of $\overline{N_k}$ is a $k$-block closure of $G$. Therefore $N_k$ is a
$k$-block closure of $G$, so that $\ell_{\mathrm{max}}(G^+)\leq \ell_{\mathrm{max}}(N_k)$. This, together with the fact that $G^+$ is a
$2$-connected spanning network on $X$ utilising $k$ Steiner points implies that $G^+$ is a $(2,k)$-MBSN. $\ \square$\end{pf}

\section{Algorithm for $k=1$}\label{sec2}
For any connected graph $G$ that is not a block, let $r(G)$ be the radius of the SCSD $C(G)$ on the set of vertices
$\displaystyle\bigcup_{Y_i\in \mathcal{Y}_0(G)} V(Y_i^*)$, where two vertices are the same colour if and only if they belong to the same
leaf-block of $G$. Let $G^{\mathrm{SD}}$ be the graph that we obtain from $G$ by introducing a Steiner point $s_0$ as follows. We locate $s_0$
at the centre of $C(G)$, and for each $Y_i\in \mathcal{Y}_0(G)$ we add an edge $s_0x$ for some $x\in Y_i^*$ where $|s_0x|\leq |s_0y|$ for all
$y\in Y_i^*$. If $G$ is $2$-connected then, to get $G^{\mathrm{SD}}$, we place $s_0$ at the midpoint of any edge $e$ of $G$ and add edges
incident to $s_0$ and the endpoints of $e$; in other words $C(G)$ will be the circle centred at the midpoint of $e$ with $r(G)=\frac{1}{2}|e|$.
Similarly to Lemma \ref{lembleaf} we have the following lemma.

\begin{lem}\label{sub1}If $G_1$ is a connected edge subgraph of $G_2$ then $r(G_1)\geq r(G_2)$.
\end{lem}

\begin{props}For any connected graph $G$, \label{isopt}$G^{\mathrm{SD}}$ is an optimal $1$-block closure of $G$.
\end{props}
\begin{pf}
This is clearly true if $G$ is a block, so assume that $G$ is connected but not a block. We first show that $G^{\mathrm{SD}}$ is 2-connected.
Let $u_1,u_2$ be any two vertices of $G^{\mathrm{SD}}$. If $u_1$ and $u_2$ are contained in the same block of $G$ then clearly there exists a
cycle in $G^{\mathrm{SD}}$ containing them both. Suppose next that $u_1$ and $u_2$ are contained in different leaf-blocks of $G$. Let
$u_i^\prime$ be a neighbour of $s_0$ in the interior of the block of $G$, say $Y_i$, containing $u_i$. We assume that $u_i,u_i^\prime,\tau(Y_i)$
are distinct, but the reasoning is similar if any of them coincide. Let $C_i$ be a cycle in $Y_i$ containing $u_i,u_i^\prime,\tau(Y_i)$. Then
there exists a path $P_i$ in $C_i$ connecting $u_i^\prime$ and $\tau(Y_i)$ and containing $u_i$. Let $P$ be a path in $G$ connecting $\tau(Y_1)$
and $\tau(Y_2)$ (note that $P$ may consist of a single vertex). Therefore the cycle formed by $P_1,P,P_2$ and the two Steiner edges incident to
$u_1^\prime$ and $u_2^\prime$ contains $u_1$ and $u_2$. The case when one of the $u_i$ coincides with $s_0$ or is contained in a non-leaf-block
is similar, and therefore for every pair of vertices of $G^{\mathrm{SD}}$ there exists a cycle containing them. Therefore $G^{\mathrm{SD}}$ is
$2$-connected.

Let $G^+$ be any optimal $1$-block closure of $G$ with Steiner point $s$. Now suppose to the contrary that $\ell_{\max}(G^+)<
\ell_{\max}(G^{\mathrm{SD}})$. Then $\ell_{\max}(G)\leq \ell_{\max}(G^+)< \ell_{\max}(G^{\mathrm{SD}})$. Then $s_0$ must be an endpoint of
$e_{\max}(G^{\mathrm{SD}})$, and therefore $\ell_{\max}(G^{\mathrm{SD}})=r(G)$. Let $C$ be the circle centred at $s_0$ and of radius
$r^\prime=\max\{|sx|\,:\,sx\mathrm{\ is\ an\ edge\ of\ }G^+\}$. Then, by Lemma \ref{lemLeaf}, $C$ is a colour-spanning disk on the interiors of
the leaf-blocks of $G$. Therefore $\ell_{\max}(G^{\mathrm{SD}})=r(G)\leq r^\prime\leq \ell_{\max}(G^+)$, which is a contradiction. $\
\square$\end{pf}

Algorithm \ref{alg1} constructs a $(2,1)$-MBSN on a set $X$ of vertices embedded in the Euclidean plane.

\begin{algorithm}[h!]
\caption{Construct a $(2,1)$-MBSN} \label{alg1}\algsetup{indent=1.5em}
\begin{algorithmic}[1]
\REQUIRE A set $X$ of $n$ vertices embedded in the Euclidean plane

\ENSURE A $(2,1)$-MBSN on $X$

\STATE Construct the $2$-RNG $R$ on $X$

\STATE Let $L$ be the ordered set of edge-lengths occurring in $R$, where ties have been broken randomly

\STATE Let $t$ be a median of $L\hspace{0.5in}$//\textit{a binary search now commences}

\REPEAT

\STATE Construct the BCF of $G_t=R(t)$

\IF {$b(G_t)> 5$ or $G_t$ is not connected}

\STATE Exit the loop and let $t$ be the median of the next larger interval

\ENDIF

\STATE Construct $C(G_t)$

\IF {$r(G_t)\leq t$}

\STATE Let $t$ be the next smaller median

\ELSE

\STATE Let $t$ be the next larger median

\ENDIF

\UNTIL no smaller value of $\max\{r(G_t),t\}$ can be found

\STATE Let $t^*\in L$ be the value that produces the minimum $\max\{r(G_t),t\}$, and let $s^*$ be the centre of $C(G_{t^*})$.

\STATE Construct a $(2,0)$-MSBN on $X\cup \{s^*\}$ and output this as the final solution
\end{algorithmic}
\end{algorithm}

\begin{thm}\label{thrm1}Algorithm \ref{alg1} correctly computes a (2,1)-MBSN on $X$ in a time of $O(n^2)$.
\end{thm}
\begin{pf}
Observe first that by Proposition \ref{isopt} for every $G_t$ Algorithm \ref{alg1} correctly computes the location of the Steiner point and the
length of the longest edge in an optimal $1$-block closure of $G_t$. Let $t_{\mathrm{opt}}=\ell_{\max}(\overline{N_1})$ and let
$G_{\mathrm{opt}}=R(t_{\mathrm{opt}})$. Note that $t_{\mathrm{opt}}\in L$, $G_{\mathrm{opt}}$ is connected since $N_1$ is connected, and (by
Corollary \ref{maxdeg}) $b(G_{\mathrm{opt}})\leq 5$. By Corollary \ref{propR}, $G_{\mathrm{opt}}^{\mathrm{SD}}$ is a (2,1)-MBSN on $X$. Any
$t\in L$ such that $G_t$ is connected, $b(G_t)\leq 5$, and $G_t^{\mathrm{SD}}$ is a $(2,1)$-MBSN on $X$ is referred to as \textit{feasible}.

Now let $t\in L$ be some value considered in the binary search. If $G_t$ is not connected or $b(G_t)>5$ then, by Lemma \ref{lembleaf}, there
exists a feasible $t^\prime$ such that $t^\prime>t$. If $r(G_t)\leq t$ then clearly there exists a feasible $t^\prime$ such that $t^\prime\leq
t$, and if $r(G_t)>t$ then, by Lemma \ref{sub1}, there exists a feasible $t^\prime$ such that $t^\prime\geq t$. Therefore  a feasible $t^\prime$
will be located by the binary search by decreasing $t$ if $r(G_t)\leq t$ and $G_t$ is connected, and increasing $t$ otherwise.

To prove the required complexity, note that the constructions of the $2$-RNG and the $(2,0)$-MBSN in Lines (1) and (15) respectively each
requires $O(n^2)$ time. The binary search in Lines (4)--(13) is on $O(n)$ elements and therefore terminates in $O(\log n)$ steps. In each step a
BCF on $G_t$ is constructed in Line (5), requiring $O(n)$ time, and an SCSD is constructed in Line (8), requiring $O(n\log n)$ time. Therefore
the total time for the search to terminate is $O(n\log^2n)$, and the total complexity is $O(n^2)$. $\ \square$\end{pf}

\section{Algorithm for $k=2$}\label{sec3}
Let $G$ be any graph on $X$ and let $G^+$ be any optimal $2$-block closure of $G$ with Steiner point set $S_2=\{s_1,s_2\}$. For any
$i\in\{1,2\}$ we denote $3-i$ by $\overline{i}$. If $G$ is a block then the construction of an optimal $2$-block closure of $G$ is easily
achieved. If $G$ is not a block but $G^+-s_i$ is a block for some $i$ (in which case $G$ is connected) then the following modification to $G^+$
will destroy this property without changing the length of the longest edge. Let $e=s_{\overline{i}}\,y$ be any Steiner edge of $G^+-s_i$. We
remove $s_i$ and edge $e$ from $G^+$, then reintroduce $s_i$ at the midpoint of line segment $s_{\overline{i}}\,y$ by adding edges $s_1s_2$ and
$s_iy$. Therefore throughout this section we assume that neither $G$ nor $G^+-s_i$ are blocks for any $i$.

\subsection{Critical edges of $G^+$}

We begin by proving a lemma that, combined with Lemma \ref{lemLeaf}, specifies a set of Steiner edges that necessarily occur in $G^+$. These
edges together with $G$ induce a subgraph of $G^+$ with a simple structure, which we then use to determine additional critical edges of $G^+$.
The benefit of knowing the critical edges becomes apparent in Section \ref{sec31}, where we present a method for locating the Steiner points of
an optimal $2$-block closure by constructing SCSDs on the blocks of $G$ containing the endpoints of the critical edges.

\begin{lem}\label{lemInd}For every isolated component $W$ of $G$ there exists a pair of Steiner $W$-edges in $G^+$. If $W$ is not a vertex
there exists a pair of independent Steiner $W$-edges in $G^+$.
\end{lem}
\begin{pf}
Clearly there exist at least two Steiner $W$-edges. Suppose that $W$ is not an isolated vertex and that no pair of independent Steiner $W$-edges
exist. Without loss of generality let $e=xs_1$ be any Steiner $W$-edge. Then either (1) all Steiner $W$-edges are incident to $x$ or (2) they
are all incident to $s_1$. If (1) is true then $x$ separates $W$ from $S_2$ in $G^+$, and if (2) is true then $s_1$ separates $W$ from $s_2$ in
$G^+$. In either case $G^+$ is not $2$-connected, which is a contradiction. Therefore an independent pair of Steiner $W$-edges must exist. $\
\square$\end{pf}

Let $E_0$ be a maximal set of external Steiner edges of $G^+$ such that: (1) every $e\in E_0$ is incident to $Y^*$ for some $Y\in
\mathcal{Y}_0(G)$ or to an isolated block of $G$, (2) no two edges of $E_0$ are incident to the same leaf-block, (3) for every isolated block
$W$ of $G$ there exists exactly two edges of $E_0$ incident to $W$ which, unless $W$ is a vertex, are independent. The set $E_0$ is referred to
as a \textit{base edge-set} for $G^+$, and its existence is guaranteed by the previous lemma and Lemma \ref{lemLeaf}. Let $E_0^{\,\prime}$ be
the set of Steiner edges not in $E_0$ and let $M_0=G^+-E_0^{\,\prime}$. If, for a given (non-block) isolated component $W$ of $G$, each edge of
$E_0$ incident to $W$ is also incident to the same Steiner point $s_i$ for some $i\in\{1,2\}$, then $W$ is called an \textit{$s_i$-covered
component}. Note that $G$ itself cannot be $s_i$-covered for some $i$ since then $G^+-s_{\overline{i}}$ would be a block. Let $M_0^\prime$ be
the subgraph of $M_0$ induced by $S_2$ and all components of $G$ that are not $s_i$-covered for any $i$.

\begin{props}\label{propcat}One of the following is true: (1) $M_0^\prime$ consists of two isolated Steiner points, (2) $M_0^\prime$ is a block,
or (3) the BCF of $M_0^\prime$ is a path with end-blocks $Y_1,Y_p$ such that $s_1\in Y_1^*$ and $s_2\in Y_p^*$.
\end{props}
\begin{pf}
If $G$ is not connected and every component is $s_i$-covered for some $i$ then clearly $M_0$ consists of exactly two isolated components and
therefore (1) holds. So let us assume that some component $W$ of $G$ is not $s_i$-covered (note that $W$ may be an isolated block). Then $s_1$
and $s_2$ are connected in $M_0^\prime$ by a path with all its internal vertices contained in $W$. Since every component of $G$ is adjacent to
at least one of the $s_i$ through an edge of $E_0$, we see that $M_0^\prime$ (and indeed $M_0$) is connected. Now suppose that $M_0^\prime$ is
not a block and that there exists a leaf-block of $M_0^\prime$, say $Y$, such that neither $s_1$ nor $s_2$ are in $Y^*$. Since $Y$ is a
leaf-block it contains at most one cut-vertex of $M_0^\prime$. If this cut-vertex is a Steiner point, say $s_1$, then $Y-s_1$ is an isolated
component of $G$ which is adjacent only to $s_1$ in $E_0$; this contradicts the definition of $M_0^\prime$. Otherwise, if $Y\cap S_2$ is empty
then $Y$ is a leaf-block of some component of $G$, and no edge in $E_0$ is incident to $Y^*$; this contradicts the choice of $E_0$. Therefore
(3) holds and the proposition follows. $\ \square$\end{pf}

\begin{cors}\label{corcat}If $G$ is connected then either $M_0$ is $2$-connected or its BCF is a path.
\end{cors}
\begin{pf}
Observe that $M_0^\prime=M_0$ in this case. $\ \square$\end{pf}

\begin{cors}\label{corcat2}If $G$ is not connected then either $G$ contains an $s_i$-covered component or $M_0$ is $2$-connected.
\end{cors}
\begin{pf}
If $G$ contains at least two components that are not $s_i$-covered then, using similar reasoning to the proof of Proposition \ref{isopt} where
it was shown that $G^{\mathrm{SD}}$ is $2$-connected, we can show that $M_0^\prime$ is $2$-connected. $\ \square$\end{pf}

\begin{figure}[htb]
  \begin{center}
    \includegraphics[scale=0.4]{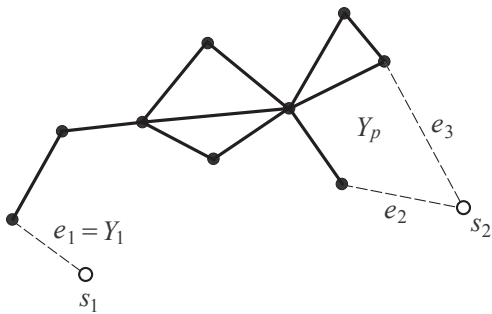}\\
  \end{center}
  \caption{$E_0=\{e_1,e_2,e_3\}$ and the BCF of $M_0=M_0^\prime$ is a path}
  \label{figPath}
\end{figure}

\begin{figure}[htb]
  \begin{center}
    \includegraphics[scale=0.4]{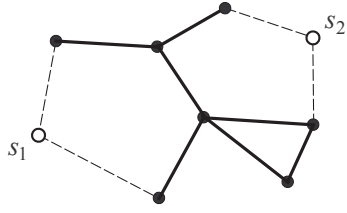}\\
  \end{center}
  \caption{$M_0=M_0^\prime$ is a block}
  \label{figBlock}
\end{figure}

\begin{figure}[htb]
  \begin{center}
    \includegraphics[scale=0.4]{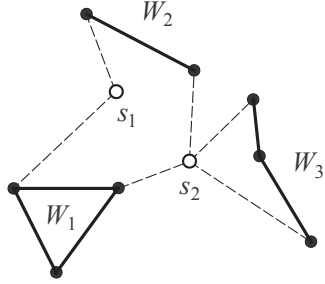}\\
  \end{center}
  \caption{$W_3$ is $s_2$-covered and $M_0^\prime$ is a block}
  \label{figsiCov}
\end{figure}

In Figs. \ref{figPath} and \ref{figBlock} we illustrate the case when $G$ is connected. Depending on the choice of $E_0$ we either attain an
$M_0$ that has a path BCF as in Fig. \ref{figPath}, or we attain an $M_0$ which is a block as in Fig. \ref{figBlock}. An example where $G$ is
not connected and contains an $s_2$-covered component is shown in Fig. \ref{figsiCov}. In this figure $G$ also contains two isolated blocks
$W_1,W_2$. In all three figures the Steiner points are represented by unfilled circles, vertices of $G$ by black filled circles, edges of $G$ by
solid lines, and edges of $E_0$ by broken lines.

As we will prove later, all critical edges of $G^+$ are specified by Lemma \ref{lemLeaf} and Lemma \ref{lemInd}, barring one particular case.
The following notation is used for this case throughout the rest of the paper. Suppose that $G$ is connected but that $M_0$ is not a block. As
per Proposition \ref{propcat} let $Y_1,...,Y_p$ be the blocks of $M_0$ as they appear in the path of the BCF, with $s_1\in Y_1^*$ and $s_2\in
Y_p^*$, and recall that $E_0^{\,\prime}$ is the set of Steiner edges of $G^+$ not contained in $E_0$. For every $i\in\{1,...,p-1\}$ let
$\tau_i=V(Y_i)\cap V(Y_{i+1})$, i.e., $\tau_i$ is the unique cut-vertex of $M_0$ common to $Y_i$ and $Y_{i+1}$. Let $B_1,...,B_p$ be the
sequence of subgraphs of $M_0$ such that $B_1=Y_1$ and for every $i\in\{2,...,p\}$, $B_i=Y_i-\tau_{i-1}$. Note then that $B_p=Y_p^*$, every
$B_i$ contains at most one cut-vertex of $M_0$, and $\{V(B_i)\}$ partitions $V(M_0)$.

\begin{lem}\label{crossedge}$E_0^{\,\prime}$ contains at least one of the following.
\begin{enumerate}
    \item An edge $s_1x$ where $x\in Y_p^*$,
    \item An edge $s_2y$ where $y\in Y_1^*$,
    \item Two edges $s_1x_1,s_2x_2$ where $x_1\notin Y_1\cup Y_p^*$; $x_2\notin Y_p\cup Y_1^*$; $x_1$ and $x_2$ are not the same cut-vertex
    of $M_0$; and if $x_1\in B_{j_1}$ and $x_2\in B_{j_2}$, then $j_2\leq j_1$.
\end{enumerate}
\end{lem}
\begin{pf}
Observe that the case when $s_1s_2$ is an edge of $E_0^{\,\prime}$ is contained in (1) or (2). Since $M_0$ is not a block $E_0^{\,\prime}$
cannot be empty. Let $i\in \{1,...,p-1\}$. Then $\tau_i$ separates $H_1=\displaystyle\bigcup_{j\geq i+1} B_j$ from
$H_2=\displaystyle\bigcup_{j\leq i} B_j-\tau_i$ in $M_0$, and therefore in $G^+$ there exists an edge connecting $H_1$ and $H_2$. Since this
edge must belong to $E_0^{\,\prime}$ (i.e., it is a Steiner edge), and there exists an edge like this for every cut-vertex of $M_0$, the result
follows. $\ \square$\end{pf}

This subsection described a number of edges (or rather, types of edges) that are necessary for a $2$-block closure of $G$. In the next
subsection we will prove that these types of edges are also sufficient.

\subsection{Constructing an optimal $2$-block closure of $G$}\label{sec31}
The construction of an optimal $1$-block closure described in Algorithm \ref{alg1} consists of locating the Steiner point at the centre of the
SCSD on the interiors of the leaf-blocks of $G$. We can also view this construction in another way. Suppose that $G$ is connected and let $M$ be
a graph topology containing $G$, a Steiner point $s$, and exactly one Steiner edge for each leaf-block of $G$. The location and precise
neighbours of $s$ in $G$ are not yet specified, yet we know that if the interior of every leaf-block of $G$ contains an endpoint of a Steiner
edge of $M$ then $M$ must be $2$-connected. Any $1$-block closure of $G$ must contain $M$, therefore by optimally \textit{embedding} $M$ (i.e.,
by determining the precise neighbours and location of $s$) we produce an optimal $1$-block closure of $G$. Our generalisation to $2$-block
closures also defines $M$ in this informal sense, but $M$ can be defined formally by, for instance, replacing each block-interior by a unique
vertex (note that block cut-vertex decompositions are often considered in this way, see \cite{timo}). Since $M$ is essentially the topology of a
graph that is obtained by removing all non-critical Steiner edges from some $2$-block closure of $G$, we refer to $M$ as a \textit{critical
topology}.

The topology of $M$ when $k=1$ can only take one general form, but when $k=2$ we will need to consider a number of candidate critical
topologies, and calculate an optimal pair of Steiner point locations for each one. The process of building a critical topology begins with the
selection of a base edge-set $E_0$. If $s_i$ is incident to $e$ in $E_0$, and $V_e$ is the block containing the other end-point of $e$, then
both $s_i$ and $e$ are said to be \textit{associated} with $V_e$. With $M_0$ defined as before we utilise Proposition \ref{propcat} to determine
whether additional Steiner edges are necessary for completing the critical topology $M$.

Once $M$ is specified, the Steiner points are located using SCSDs and \textit{farthest colour Voronoi diagrams} (FCVDs). The FCVD is defined in
\cite{abel} as follows. Let $\mathcal{C}=\{P_1,...,P_q\}$ be a collection of $q$ sets of $n$ coloured points. If $p\in P_i$, i.e., $p$ is a
point of colour $i$, we put all points of the plane in the \textit{region} of $p$ for which $i$ is the farthest colour, and $p$ the nearest
$i$-coloured point. In other words, $z$ belongs to the region of $p$ if and only if the closed circle centred at $z$ that passes through $p$
contains at least one point of each colour, but no point of colour $i$ is contained in its interior. The FCVD for $\mathcal{C}$ is the
decomposition of the plane into these regions; in other words the edges and vertices of the FCVD are the intersections of boundaries of regions.

\begin{thm}(see \cite{abel}) For constant $q$ an FCVD on $\mathcal{C}$ can be computed in $O(n^2)$ time, and its structural complexity is $O(n)$.
\end{thm}

\begin{cors}(see \cite{abel}) Given the FCVD, an SCSD on $\mathcal{C}$ can be found in $O(n)$ time.
\end{cors}
\begin{pf}The centre of the SCSD is either a vertex or the midpoint of an edge of the FCVD.
$\ \square$\end{pf}

Let $C$ be an SCSD on $\mathcal{C}$ and let $x$ be the centre of $C$.

\begin{lem}A set $D(x)$ of cardinality $q$ containing a closest point of each colour to $x$ can be constructed in
$O(n\log n)$ time.
\end{lem}
\begin{pf}
A closest point of $P_i$ is found by constructing a standard Voronoi diagram on $P_i$ and then performing point-location on $x$. $\
\square$\end{pf}

Due to the previous result we assume in the rest of this section that the set $D(x)$ is known after any construction of an SCSD. It will be seen
later that the purpose of $D(x)$ is to specify the neighbours of the Steiner points.

Recall that we are assuming that $G$ is not a block. In order to choose a candidate base edge-set $E_0$ we partition the set $\mathcal{Y}_0(G)$
into two sets $\mathcal{P}=\{\mathcal{Y}^1,\mathcal{Y}^2\}$, where one of the sets may be empty if $G$ is not connected. Let $\mathcal{Z}$ be
the set of isolated blocks of $G$. In $E_0$ we then associate $s_1$ with each member of $\mathcal{Y}^1$, and $s_2$ with each member of
$\mathcal{Y}^2$. Each $s_i$ is also associated with every member of $\mathcal{Z}$. The edge-set $E_0$ defines the graph $M_0$ (as in the
previous subsection) We now discuss three different cases depending on the structure and connectivity of $M_0$. In each case we show how to
construct a critical topology $M$ and how to embed $M$ optimally.\\

\noindent\textbf{Case 1}: $M_0$ is $2$-connected.\\
In this case no additional edges are required for an optimal $2$-block closure of $G$, therefore we let $M=M_0$. Suppose first that
$|\mathcal{Z}|=0$. We assign a unique colour to each $Y^*$ where $Y\in\mathcal{Y}^1$. Let $s_1$ be the centre of the SCSD on these colour sets.
We then perform a similar operation in order to find the location of $s_2$. When $|\mathcal{Z}|\neq 0$ we need to make sure that $V(Z)\cap
D(s_1)\cap D(s_2)=\emptyset$ for every $Z\in \mathcal{Z}$ with $|V(Z)|>1$. This is because $D(s_i)$ specifies the neighbours of $s_i$ in the
optimal embedded version of $M$, and, by the choice of $E_0$, if $Z$ is not a vertex then $s_1$ and $s_2$ must have distinct neighbours in $Z$.
If $Z$ is an isolated vertex then it will be assigned a unique colour along with the leaf-blocks of $\mathcal{Y}^i$ when locating each $s_i$,
therefore for the remainder of Case 1 we assume that none of the members of $\mathcal{Z}$ are vertices.

Next suppose that $|\mathcal{Z}|=1$. We proceed exactly as before in order to locate $s_1$. Let $Z\in \mathcal{Z}$ and let $y= V(Z)\cap D(s_1)$.
When locating $s_2$ we proceed as before, but this time we do not include $y$ when colouring $Z$. Next the entire process is repeated, but this
time $s_2$ is located before $s_1$. The cheapest of these two solutions (determined by the largest radius of the two SCSDs) is picked as the
final solution.

The final subcase we consider is when $|\mathcal{Z}|=5$, so that each $\mathcal{Y}^i$ is empty. Our method is essentially a generalisation of
the previous subcase, and all other subcases are subsumed by it. Suppose that $D(x)=\{y_i\in Z_i\}$, where $\mathcal{Z}=\{Z_i\}$ and $x$ is the
centre of the SCSD on $\mathcal{Z}$.\\

\noindent\textbf{Claim}: For some $i\in\{1,2\}$ there exists an SCSD $C_i$ such that the optimal location of $s_i$ is the centre of $C_i$, and
such that at least one member of $D(s_i)$ is contained in $\{y_i\}$. By symmetry we may assume that $i=1$.
\begin{pf}
If this were not true then we could relocate $s_2$ at $x$, and let the neighbour-set of $s_2$ be $\{y_i\}$ in the embedded version of $M$.
Clearly this will not increase the length of any edge and $V(Z_i)\cap D(s_1)\cap D(s_2)$ will be empty for every $Z_i$. $\ \square$\end{pf}

For every $j\in\{1,...,5\}$ we perform the following process. Suppose without loss of generality that $j=1$. Let $C_1^\prime$ be the SCSD, with
centre $x_1$, on $\{y_1\},Z_2,..,Z_5$ and let $C_2^\prime$ be the SCSD, with centre $x_2$, on $Z_1-\{y_1\},Z_2,...,Z_5$. Similarly to the
previous claim, we may assume that $D(s_i)\cap D(x_{j_1})\cap V(Z_{j_2})\neq \emptyset$ for some $i,j_1\in\{1,2\}$, and some $j_2\in
\{2,...,5\}$, where $s_i$ is an optimal Steiner point location. We perform the following process for every such $j_1,j_2$ and $y^\prime\in
D(x_{j_1})\cap V(Z_{j_2})$. Suppose without loss of generality that $y^\prime\in D(x_1)\cap V(Z_2)$. Let $C_1^{\prime\prime}$ be the SCSD on
$\{y_1\},\{y^\prime\},Z_3,...,Z_5$ and let $C_2^{\prime\prime}$ be the SCSD on $Z_1-\{y_1\},Z_2-\{y^\prime\},Z_3,...,Z_5$, and continue the
process as before. The process ends when we have located $s_1$ and $s_2$ such that $D(s_1)\cap D(s_2)=\emptyset$. The optimal embedded version
of $M$ is selected as a cheapest solution of all the various iterations.
The total time-complexity in Case 1 is $O(n\log n)$.\\

\noindent\textbf{Case 2}: $M_0$ is not $2$-connected and there are no $s_j$-covered components of $G$ for any $j\in\{1,2\}$.\\
By Corollary \ref{corcat2} this case only arises when $G$ is connected. There are two subcases here, and we consider both before picking
a cheapest solution.\\

\noindent\textbf{Subcase 2.1}: Edge $s_1s_2$ is not included in $M$.\\
We use the notation from Lemma \ref{crossedge}. If $Y_1$ consists of a single edge then let $J_1=1$, else let $J_1=\emptyset$; similarly if
$Y_p$ consists of a single edge then let $J_2=p$, else let $J_2=\emptyset$. Let $i\in \{1,...,p\}-J_1-J_2$. If $i=p$ then let $E_0^{\,\prime}$
consist of a single edge incident to $s_1$ and associated with $Y_p^*-s_2$. If $i=1$ then let $E_0^{\,\prime}$ consist of a single edge incident
to $s_2$ and associated with $Y_1^*-s_1$. Otherwise, let $E_0^{\,\prime}$ consist of two edges $e_1,e_2$, where $e_1$ is incident to $s_1$ and
associated with $B_i$, and $e_2$ is incident to $s_2$ and associated with $\displaystyle\bigcup_{j\leq i}B_j-\tau_i$.

\begin{lem}Critical topology $M=M_0+E_0^{\,\prime}$ is $2$-connected for any $i\in \{1,...,p\}-J_1-J_2$.
\end{lem}
\begin{pf}
Clearly $M$ is connected. Since $M_0$ is a connected edge-subgraph of $M$, if $x$ is a cut-vertex of $M$ then $x$ is also a cut-vertex of $M_0$.
Therefore, if $x$ is a cut-vertex of $M$ then $x=\tau_j$ for some $j\in\{1,...,p-1\}$, so that $x$ separates $H_1=\displaystyle\bigcup_{j_0\geq
j+1} B_{j_0}$ from $H_2=\displaystyle\bigcup_{j_0\leq j} B_{j_0}-x$ in $M$. But by the definition of $E_0^{\,\prime}$ either $i\geq j+1$ and
$e_1\in E_0^{\,\prime}$ is associated with $B_i$, or $i\leq j$ and $e_2\in E_0^{\,\prime}$ is associated with $\displaystyle\bigcup_{j_0\leq
i}B_{j_0}-\tau_i$. In either case there is an edge of $E_0^{\,\prime}$ connecting a vertex of $H_1$ and a vertex of $H_2$. Therefore no such
separating vertex $x$ exists. $\ \square$\end{pf}

For locating the Steiner points we assume that $|E_0^{\,\prime}|=2$, the other case is similar. Let $I_0=\{1,...,p\}-J_1-J_2$. We perform a
binary search on $I_0$ in order to find the cheapest solution of the following form. Let $a\in I_0$, let $H_1^a=\displaystyle\bigcup_{j\geq
a}B_j$, and let $s_1$ be located at the centre of the SCSD on the members of $\mathcal{Y}^1$ and on $H_1^a$. To locate $s_2$ suppose that
$D(s_1)\cap H_1^a$ lies in $B_b$, where $b=b(a)\geq a$. Let $H_2^b=\displaystyle\bigcup_{j\leq b}B_j-\tau_b$ and locate $s_2$ at the centre of
the SCSD on the members of $\mathcal{Y}^2$ and on $H_2^b$. For $i=1,2$ let $r_a^i$ be the radius of the SCSD constructed for $s_i$. The binary
search on $I_0$ will find the value of $a$ for which $r_a=\max\{r_a^1,r_a^2\}$ is a minimum. Observe that there must exist an $a\in I_0$ such
that the Steiner point locations constructed by this method for $a$ are optimal for a $2$-block closure of the current type.

We begin the search with a median value of $I_0$. Suppose that the current iteration of the search is $a\in I_0$. If $r_a^1\geq r_a^2$ then we
decrease $a$ for the next iteration, otherwise we increase $a$. We repeat this until no smaller value of $r_a$ is found. To see why the search
will terminate at an optimal value of $a$ suppose first that $r_a^1\geq r_a^2$ at some iteration. Now let $a^\prime\in I_0$ such that
$a^\prime\geq a$. Then since $H_1^{a^\prime}\subseteq H_1^a$ we must have $r_{a^\prime}^1\geq r_a^1\geq r_a$. Therefore $a^0\leq a$ for some
optimal $a^0$. Next suppose that $r_a^1< r_a^2$. Then, by similar reasoning for $H_2^b$, $b(a^0)\geq b(a)$ for some optimal $a_0$. But $b$ is a
non-decreasing function of $a$, and therefore we may assume that $a^0 \geq a$.

Since $|I_0|\in O(n)$ the search will terminate in $O(\log n)$ steps. At each step we construct two SCDS, and therefore the total time to locate
the optimal Steiner point pair is $O(n\log^2n)$.\\

\noindent\textbf{Subcase 2.2}: Edge $s_1s_2$ is included in $M$.\\
Similarly to the previous subcase we have the following result:

\begin{lem}Critical topology $M=M_0+s_1s_2$ is $2$-connected.
\end{lem}

When embedding $M$ there are a few possibilities depending on the locations and the number of determinators of the SCSDs for each Steiner point,
but these cases are all similar to the results of \cite{bae2} and will therefore not be discussed in much detail.

We briefly look at one of the cases. When each Steiner point is a determinator of the other Steiner point's SCSD and both SCSDs have three
determinators, we may locate the Steiner points by constructing two FCVDs, one on the leaf-blocks in $\mathcal{Y}^1$ and another on the
leaf-blocks in $\mathcal{Y}^2$. We then select an edge of each FCVD before solving a quartic equation to locate the Steiner points. This is
possible since each of the two edges contains one of the Steiner points, and the distance between the Steiner points is equal to the common
radius of the SCSDs. The maximum time for locating two adjacent Steiner points is $O(n^2)$ since we need to consider every pair of
$O(n)$ edges.\\

\noindent\textbf{Case 3}: $M_0$ is not $2$-connected and $G$ contains at least one $s_i$-covered component for some $i\in\{1,2\}$.\\
This case only occurs when $G$ is not connected. For $j=1,2$ and a set of integer indices $I_j$ let $\{W_i^j:i\in I_j\}$ be the set of
$s_j$-covered components of $G$. Let $E_j$ be the set of edges containing exactly one edge $e_i$ for each $i\in I_j$ such that $e_i$ is incident
to $s_{\overline{j}}$ and is associated with $W_i^j$. Observe by Lemma \ref{lemInd} that $E_1$ and $E_2$ are necessarily in a $2$-block closure
of $G$.

\begin{lem}Critical topology $M=M_0+E_1+E_2$ is $2$-connected.
\end{lem}
\begin{pf}
Observe that $M$ is connected since the addition of any edge of $E_1$ or $E_2$ to $M_0$ creates a path connecting $s_1$ and $s_2$. Suppose to
the contrary that $M$ has a cut-vertex $x$. Then $x$ is also a cut-vertex of $M_0$ and is therefore one of the following vertices: (1) a
cut-vertex of $M_0^\prime$, (2) a Steiner point, (3) a non-Steiner end-point of a Steiner $V$-edge in $E_0$, where $V$ is an $s_i$-covered
component. Suppose that (1) holds and suppose without loss of generality that $W$ is an $s_2$-covered component of $G$. Note that $x$ separates
$s_1$ and $s_2$ in $M_0$, and therefore also separates these vertices in $M$. Let $e\in E_2$ be a Steiner $W$-edge incident to $s_1$, and let
$e^\prime\in E_0$ be a Steiner $W$-edge incident to $s_2$. Let $P_1$ be a path in $W$ connecting the non-Steiner end-points of $e$ and
$e^\prime$, and let $P_2$ be a path in $M_0$ connecting $s_1$ and $s_2$ (and therefore containing $x$). Then $P_1,P_2$ and the edges
$e,e^\prime$ form a cycle in $M$ containing $s_1,s_2$ and $x$, which contradicts the fact that $x$ separates $s_1$ and $s_2$. Cases (2) and (3)
are handled similarly since in these cases the cut-vertices lie on the same type of cycle. Therefore the lemma follows. $\ \square$\end{pf}

To find the location of $s_i$ we assign a unique colour to every $W_j^{\overline i}$ and to each $Y\in \mathcal{Y}^i$ and $Z\in \mathcal{Z}$. We
then proceed similarly to Case 1, and again consider subcases depending on the cardinality of $|\mathcal{Z}|$. The sets $W_j^i$ are treated
exactly as leaf-blocks are in Case 1. The total run-time is therefore also $O(n\log n)$.

The above three cases cover all possibilities. To close this section we observe that the pair of Steiner point locations $S_2=\{s_1,s_2\}$
produced in the relevant case will be optimal for the embedded version of $M$. In other words, for any optimal $2$-block closure $G^+$ of $G$
such that $G^+$ contains the critical topology $M$ (and note that we have shown it must contain $M$ for one of the cases), the embedded version
of $M$ is an optimal $2$-block closure of $G$. The proof of this fact is similar to the second part of the proof of Proposition \ref{isopt}, and
we therefore do not provide further details.

For any given $G$ and some $M$ let $r(M)$ be the maximum radius of an SCSD used to optimally embed $M$. Let $r(G)=\min\{r(M)\}$ and let
$G^{\mathrm{SD2}}$ be an optimally embedded $M$ attaining $r(G)$. Then clearly $G^{\mathrm{SD2}}$ is an optimal $2$-block closure of $G$.
Similarly to Lemma \ref{sub1} we have the following result.

\begin{lem}\label{sub2}If $G_1$ is an edge subgraph of $G_2$ then $r(G_1)\geq r(G_2)$.
\end{lem}

We present Algorithm \ref{alg2} for constructing a $(2,2)$-MBSN.

\begin{algorithm}[h!]
\caption{Construct a $(2,2)$-MBSN} \label{alg2}\algsetup{indent=1.5em}
\begin{algorithmic}[1]
\REQUIRE A set $X$ of $n$ vertices embedded in the Euclidean plane

\ENSURE A $(2,2)$-MBSN on $X$

\STATE Construct the $2$-RNG $R$ on $X$

\STATE Let $L$ be the ordered set of edge-lengths occurring in $R$, where ties have been broken randomly

\STATE Let $t$ be a median of $L$

\REPEAT

\STATE Construct the BCF of $G_t=R(t)$

\IF {$b(G_t)> 10$}

\STATE Exit the loop and let $t$ be the median of the next larger interval

\ENDIF

\FORALL{valid partitions $\mathcal{P}=\{\mathcal{Y}^1,\mathcal{Y}^2\}$ of $\mathcal{Y}_0(G_t)$}

\STATE Let $E_0$ be the base edge-set determined by $\mathcal{P}$ and the isolated blocks of $G_t$

\STATE Construct the BCF of $M_0$

\STATE Use the structure of the BCF of $M_0$ to determine the critical topology $M$ and its optimal embedding, by calling the relevant procedure
from Case 1 -- 3

\ENDFOR

\IF {$r(G_t)\leq t$}

\STATE Let $t$ be the next smaller median

\ELSE

\STATE Let $t$ be the next larger median

\ENDIF

\UNTIL no smaller value of $\max\{r(G_t),t\}$ can be found

\STATE Output the embedded $M$ producing the minimum $\max\{r(G_t),t\}$
\end{algorithmic}
\end{algorithm}

\newpage
\begin{thm}Algorithm \ref{alg2} correctly computes a $(2,2)$-MBSN on $X$ in a time of $O(n^2\log n)$.
\end{thm}
\begin{pf}
The correctness proof is similar to that of Theorem \ref{thrm1}. Let $t_{\mathrm{opt}}=\ell_{\max}(\overline{N_2})$ and
$G_{\mathrm{opt}}=R(t_{\mathrm{opt}})$. Then $G_{\mathrm{opt}}^{\mathrm{SD2}}$ is a $(2,2)$-MBSN on $X$ and we proceed as before.

To prove complexity we note that the longest time that arises during the binary search is $O(n^2)$ in Line (11), Subcase 2.2 when the Steiner
points are adjacent to each other. Iterating through all valid partitions in Line (8) requires constant time, and constructing the BCF of $M_0$
in Line (10) takes at most $O(n)$ time $\ \square$\end{pf}

It should be noted that it is possible to replace all occurrences of the $2$-RNG in Algorithm \ref{alg2} with the complete graph on $X$, without
altering the essential nature of the algorithm. Since each iteration of the algorithm already requires $O(n^2)$ time, and the main difference in
complexity in the two versions is the time required to produce the BCFs, the final complexity would still be $O(n^2\log n)$. Even though the
limiting complexity remains unchanged, using the complete graph will become an issue during practical implementations because the BCF is
constructed so often. For this reason, and for the sake of symmetry with the $k=1$ case, we make use of the $2$-RNG here.

\section{Conclusion}\label{sec4}
By using properties of $2$-connected graphs, $2$-relative neighbourhood graphs, and smallest colour spanning disks, we produced two fast and
exact polynomial time algorithms for solving the Euclidean bottleneck $2$-connected $k$-Steiner network problem when $k=1,2$. Fundamental to our
algorithms is the fact that any graph can be uniquely decomposed into blocks such that the resulting graph is a forest. This allowed us to
characterise the set of edges which occur in an optimal solution. The properties of these edges are crucial in determining the colour sets upon
which the spanning disks should be constructed. In turn, the spanning disks determine the locations of the optimal Steiner points. In the $k=1$
case this gave us an algorithm of complexity $O(n^2)$, and $O(n^2\log n)$ when $k=2$.

Regarding the $k\leq 2$ problem on other planar norms, observe that our connectivity related results are based on topological properties, and
therefore hold for all metrics. Smallest colour-spanning disks and farthest colour Voronoi diagrams find analogs the $L_p$ planes: see
\cite{abel,hutt}. A generalisation of the $2$-relative neighbourhood graph to $L_p$ norms has not been considered in the literature, however
algorithms do exist for the construction of $1$-relative neighbourhood graphs in these planes (see \cite{jar}). It might be possible to extend
the results of \cite{jar} but, irrespectively, replacing all occurrences of the $2$-RNG in our algorithms by the complete graph on $X$ leads to
an increase in complexity of only a $\log n$ factor when $k=1$, and no increase when $k=2$.

A future goal is to extend our results to general values of $k$ and also to graphs of higher connectivity. We believe that this can be achieved
through more sophisticated methods based on the ones developed in this paper; this is one of our current topics of research.

\end{document}